\documentclass[11pt]{article}
\usepackage{newlfont,amsfonts,amssymb,
amsmath,amsthm,amsgen,amscd,wasysym}
\usepackage{rotating}
\newcommand{\C}{\ensuremath{\mathbb{C}}}

\newcommand{\R}{\ensuremath{\mathbb{R}}}

\newcommand{\bcase}{\begin{case}}
 \newcommand{\ecase}{\end{case}}

\newcommand{\bclaim}{\begin{claim}}
\newcommand{\eclaim}{\end{claim}}

\newcommand{\bstep}{\begin{step}}
\newcommand{\estep}{\end{step}}

\newcommand{\bhlem}{\begin{hlem}}
\newcommand{\ehlem}{\end{hlem}}

\newcommand{\bleer}{\begin{leer}}
\newcommand{\eleer}{\end{leer}}
\newcommand{\bde}{\begin{de}}
\newcommand{\ede}{\end{de}}

\newcommand{\ol}{\overline}

\newcommand{\mf}{\mathfrak}
\newcommand{\bs}{\begin{satz}}
\newcommand{\es}{\end{satz}}
\newcommand{\btheo}{\begin{theo}}
\newcommand{\etheo}{\end{theo}}
\newcommand{\bfolg}{\begin{folg}}
\newcommand{\efolg}{\end{folg}}
\newcommand{\blem}{\begin{lem}}
\newcommand{\elem}{\end{lem}}
\newcommand{\bnote}{\begin{note}}
\newcommand{\enote}{\end{note}}
\newcommand{\bprf}{\begin{proof}}
\newcommand{\eprf}{\end{proof}}
\newcommand{\bd}{\begin{displaymath}}
\newcommand{\ed}{\end{displaymath}}
\newcommand{\be}{\begin{eqnarray*}}
\newcommand{\ee}{\end{eqnarray*}}
\newcommand{\eeqa}{\end{eqnarray}}
\newcommand{\beqa}{\begin{eqnarray}}
\newcommand{\bi}{\begin{itemize}}
\newcommand{\ei}{\end{itemize}}
\newcommand{\bnum}{\begin{enumerate}}
\newcommand{\enum}{\end{enumerate}}
\newcommand{\la}{\langle}
\newcommand{\ra}{\rangle}

\newcommand{\beq}{\begin{equation}}
\newcommand{\eeq}{\end{equation}}

\newcommand{\rr}{\mathbb{R}}

\newcommand{\ccc}{\mathbb{C}}

\newcommand{\turnsimeq}{\mbox{\begin{turn}{270}$\simeq$  \end{turn}}}

\newcommand{\earr}{\end{array}\]}
\newcommand{\barr}{\[\begin{array}}
\newcommand{\bvec}{\left(\begin{array}{c}}
\newcommand{\evec}{\end{array}\right)}

\newcommand{\lag}{\mathfrak{g}}
\newcommand{\lan}{\mathfrak{n}}
\newcommand{\lah}{\mathfrak{h}}

\newcommand{\+}{\oplus}

\newcommand{\rrn}{\mathbb{R}^n}
\newcommand{\laso}{\mathfrak{so}}

\newcommand{\lasl}{\mathfrak{sl}}

\newcommand{\w}{\omega}

\newcommand{\s}{\sigma}

\newcommand{\bbem}{\begin{bem}}
\newcommand{\ebem}{\end{bem}}
\newcommand{\bbez}{\begin{bez}}
\newcommand{\ebez}{\end{bez}}
\newcommand{\bbsp}{\begin{bsp}}
\newcommand{\ebsp}{\end{bsp}}

\newcommand{\lam}{\mathfrak{m}}
\newcommand{\lap}{\mathfrak{p}}

\newcommand{\lak}{\mathfrak{k}}

\newcommand{\W}{\Omega}

\newcommand{\wt}{\widetilde}

\newcommand{\e}{\mathrm{e}}
\renewcommand{\i}{\mathrm{i}}
\renewcommand{\span}{\mathrm{span}}

\renewcommand{\i}{\mathrm{i}}

\theoremstyle{definition}
\newtheorem{de}{Definition}
\newtheorem{bem}[de]{Remark}
\newtheorem{bez}[de]{Notation}
\newtheorem{bsp}[de]{Example}
\theoremstyle{plain}
\newtheorem{lem}{Lemma}
\newtheorem{satz}{Proposition}
\newtheorem{folg}{Corollary}
\newtheorem{theo}{Theorem}

\begin{document}
\newcommand{\K}{\ensuremath{\mathbb{K}}}
\bibliographystyle{alpha}
\title{Connected  subgroups of $SO(2,n)$ acting irreducibly on  $\R^{2,n}$}
\author{Antonio J. Di Scala\thanks{{\sc  Dipartimento di Matematica,
Politecnico di Torino,
Corso Duca degli Abruzzi 24, 10129 Torino, Italy}, {\tt
 antonio.discala@polito.it }}\ \  and  Thomas Leistner\thanks{{\sc  Department Mathematik, Universit\"at Hamburg, Bundesstra{\ss}e 55, D-20146 Hamburg, Germany}, {\tt leistner@math.uni-hamburg.de} \newline {\em  Date: June 16, 2008}}}
 \date{ }
\maketitle
\begin{abstract}We classify all connected subgroups of $SO(2,n)$ that act irreducibly on $\rr^{2,n}$. Apart from $SO_0(2,n)$ itself these are $U(1,n/2)$, $SU(1,n/2)$,  if $n$ even, $S^1\cdot SO(1,n/2)$ if $n$ even and $n\ge 2$, and $SO_0(1,2)$  for $n=3$. Our proof is based on the Karpelevich Theorem and uses the classification of totally geodesic submanifolds of complex hyperbolic space and of the
Lie ball. As an application we obtain a list of possible irreducible holonomy groups of Lorentzian conformal structures, namely $SO_0(2,n)$, $SU(1,n)$, and $SO_0(1,2)$. \\[1mm] 
{\em MSC:}  22E46; 
 53C35;
53C40; 
53C29;
53A30.
\\[1mm]
{\em Keywords:} Irreducible, orthogonal representations; Lie ball; complex hyperbolic space; totally geodesic submanifolds of symmetric spaces; conformal holonomy; Berger list.
 \end{abstract}

\section{Background, result, and applications}
One of the results  at the origins of modern differential geometry is Marcel Berger's classification of irreducible connected  holonomy groups of complete semi-Riemannian manifolds. The most striking feature of this {\em Berger list} is that it is rather short. This is more surprising in some signatures and more natural in others. For example,  that the only possible connected irreducible holonomy group of Lorentzian manifolds is $SO_0(1,n)$ is due to the fact that there are no proper connected subgroups of $SO_0(1,n)$ that act irreducibly on $\rr^{1,n}$ \cite{olmos-discala01}. In Riemannian signature, on the other hand, it is more surprising that only so few groups occur as holonomy groups, taking into account that any representation of a compact group, in particular irreducible ones, is orthogonal with respect to a positive definite scalar product. For a recent proof of Berger's theorem for Riemannian manifolds see \cite{Olmos05}.

In this paper we consider the case of signature $(2,n)$ and study connected subgroups of $SO(2,n)$ that act irreducibly on $\rr^{2,n}$. We give a classification  of these groups:
\btheo\label{list}
Let $G\subset SO(2,n)$ be a connected Lie group that acts irreducibly on $\rr^{2,n}$. Then $G$ is conjugated to one of the following,
\bnum
\item for arbitrary $n\ge 1$:  $SO_0(2,n)$,
\item
for $n=2p$ even:  $U(1,p)$, $SU(1,p)$, or $S^1\cdot SO_0(1,p)$  if $ p>1$,
\item for $n=3$: $SO_0(1,2) \subset SO(2,3)$, .
\enum
\etheo
Our interest in signature $(2,n)$ is twofold. One aspect is the more general interest in the Berger list. Our result shows that there is only one group, namely $S^1\cdot SO(1,n)$,  that does not appear in the Berger list, i.e. that is not a holonomy group for a metric of signature $(2,n)$.

 More important is the relation to conformal  Lorentzian structures.   To a Lorentzian conformal structure in dimension $n$, which is defined as an equivalence class Lorentzian metrics differing by a scaling function, one may assume a conformally invariant Cartan connection, the holonomy group of which is contained in $SO(2,n)$. For the so-called {\em conformal holonomy} the algebraic restrictions are much more difficult to handle than the Berger criterion in case of metric holonomy algebras. Hence, it was natural to ask first: What are possible connected subgroups of $SO(2,n)$ that act irreducibly on $\rr^{2,n}$? Our answer to this question gives a list of possible candidates for {\em special conformal Lorentzian structures}, a name which refers to --- in analogy to special Riemannian structures ---  Lorentzian conformal structures with {\em irreducibly} acting conformal holonomy group (for indecomposable, non-irreducible Lorentzian conformal structures we refer the reader to \cite{leistner05a}). Now,  two of the groups Theorem \ref{list} are known to be Lorentzian conformal holonomy groups,  $SO_0(2,n)$ itself and $SU(1,n/2)$, the first being the generic conformal holonomy, the second being that of a Fefferman space (see for example \cite{baum07}).
In \cite{leitner08} it is proven that if a {\em connected} conformal holonomy group is contained $U(1,n/2)$ then it is already contained in $SU(1,n/2)$. Hence,  $S^1\cdot SO_0(1,n/2)$ cannot occur as connected conformal holonomy group of a Lorentzian conformal structure, because it is not contained in $SU(1,n/2)$. We get the following consequence.
\bfolg
Let $G\subset SO(2,n)$ be the connected conformal holonomy group of a Lorentzian conformal structure. If $G$ acts irreducibly on $\rr^{2,n}$, then
\[ G=SO_0(2,n)\text{, or }G=SU(1,n/2)\text{ if $n$ is even, or }G=SO_0(1,2)\text{ if   $n=3$.}\]
\efolg
Unfortunately, we cannot yet exclude the exceptional case of $SO_0(1,2)\subset SO(2,3)$ as a possible conformal holonomy of a $3$-dimensional Lorentzian manifold. We only know that
$SO(1,2) $ does not define a {\em conformal Cartan reduction} in the sense of \cite[Section 3.3]{alt08}.
Such a conformal Cartan reduction of $SO(p+1,q+1)$ to a group $G\subset SO(p+1,q+1)$ exists if and only if $G$ acts transitively on the pseudo-sphere $S^{p,q}= SO(p+1,q+1)/P$, where $P$ is the parabolic subgroup defined as the stabiliser of a light-like line in $\rr^{p+1,q+1}$. Examples of conformal Cartan reductions are given by $SU(p+1,q+1)\subset SO(2p+2,2q+2)$, see \cite{baum07} or \cite{cap-gover06},  the non-compact $G_{2(2)}\subset SO(3,4)$ in \cite{nurowski04,nurowski07}, and $Spin(3,4)\subset SO(4,4)$ in \cite{bryant06}, and they are linked to so-called {\em Fefferman constructions}.
Now, the action of $SO_0(1,2)$ on $S^{1,2}=SO(2,3)/P$  is not  transitive (see Appendix \ref{so12}) and hence does not define a conformal Cartan reduction, but we do not know
if this already excludes $SO_0(1,2)\subset SO(2,3)$ as an irreducible conformal holonomy. To clarify this question lies beyond the scope of this paper and  will be subject to further studies.
%


\bigskip

Our proof of Theorem \ref{list} is based on the Theorem of Karpelevich and Mostow.
\begin{theo}\label{Karp}(Karpelevich \cite{karpelevich53}, Mostow \cite{mostow55}, also \cite{discala-olmos07})
Let $M =  {Iso}(M)/K$ be a Riemannian symmetric space of
non-compact type. Then any connected and semisimple subgroup  $G$
of the full isometry group ${Iso}(M)$ has a totally geodesic orbit
$G\cdot p \subset M$.
\end{theo}
We will apply this theorem to a connected subgroup $G$ of $SO(2,n)$ that acts irreducibly on $\rr^{2,n}$ and to  the Riemannian symmetric spaces that are related to $SO(2,n)$: the complex hyperbolic space $\C H^{n}=SU(1,n)/U(n)$ and the Grassmannian of   negative definite planes in $\rr^{2,n}$ given as $SO_0(2,n)SO(2)\cdot SO(n)$ and as  $SO(2,n)/SO(2)\cdot SO(n)$ if one considers {\em oriented} negative planes. The latter has two connected components and can be realised in $\C P^{n+1}$ as the  submanifold of  negative definite lines in $\C^{2,n}$. Its connected component is called {\em Lie ball}. In applying Karpelevich's Theorem we have to deal with two difficulties that are related to each other: First, we cannot assume that $G$ is semisimple, and secondly, if $T$ is a totally geodesic orbit with isometry group $H=Iso (T)$ our group $G$ in question can be the product of $H$ with the group $I(T)$ that is defined as
\[I(T):=\{A \in G\mid A|_T=\mathrm{Id}_T\}.\]
$I(T)$ is a normal subgroup in $Iso(T)$.
We know that $G$ is reductive but it may be that its semisimple part does not act irreducibly. On the other hand, it might happen that $T$ is the orbit of a group $H$ that does not act irreducibly but that $I(T)\cdot H$ acts irreducibly.
Overcoming these difficulties,
our proof will consist of  three main steps:
\bnum
\item Show that if $G\subset SO_0(2,n)$ acts irreducibly, then it is {\em simple} or contained in $U(1,n)$.
\item Classify connected subgroups of $U(1,n)$ acting irreducibly on $\rr^{2,n}$ using:
\bnum
\item $G$ is reductive with possible centre $S^1$,
\item By Karpelevich's Theorem
 applied to $\C H^n$, the orbits $T$ of the semisimple part are isometric to either $\C H^k$  or to real hyperbolic spaces $\rr H^k$ for $k\le n$.
\item $I(T)$ can be calculated.
\enum
\item If $G$ is not in $U(1,n)$, $G$  is simple and we apply Karpelevich's Theorem to the Lie ball
$SO_0(2,n)/SO(2)\cdot SO(n)$. Then we use
the classification of totally geodesic orbits in the complex quadric $SO(n+2)/SO(2)\cdot SO(n)$ by \cite{klein08}, transfer it by duality to the Lie ball  and obtain $G$ as isometry group of these orbits. As $G$ is simple, $I(T)$ can be ignored.
 \enum

\section{Algebraic preliminaries}
\subsection{Irreducible representations of real Lie algebras}
Most of the groups appearing in the theorem are well known. However, regarding $S^1\cdot  SO_0(1,p)$  we will make some remarks about irreducible representations of real Lie algebras and about symmetric space which will be useful in what follows. The transition to Lie algebras is justified by the restriction to connected subgroups of $SO(2,n)$.

Let $\lag$ be a real Lie algebra and $E$ an irreducible {\em real} representation. We say that $E$ is of {\em real type} if $E^\C:=E\otimes\C$ is irreducible as well.  Otherwise we say that $E$ is {\em not of real type}.  In the latter case there is a splitting of $E^\C$ as
$E^\C=V\+ \ol{V}$ where $V$ is an irreducible complex representation and $\ol{V}$ the conjugate representation w.r.t. the real form $E\subset E^\C$.
By the conjugate representation we mean the representation of $\lag$ on $\ol{V}$ defined by
\[ A\cdot \ol{v}:= \ol{A\cdot v}.\]
In deed, if $V$ is a complex invariant subspace  of $E^\C$, the complex
subspaces  $V +\ol{V}$ and $V\cap\ol{V}$ are invariant as well. On the other hand, they are equal to their complex conjugate, and thus, complexifications of real invariant subspaces. As $E$ is irreducible, we obtain that $V\+ \ol{V}=E^\C$.

On the other hand, it is $\left(V_\rr\right)^\C=V\+\ol{V}$ and multiplication with $\mathrm{i}$ defines an invariant isomorphism $J$ on $V_\rr$ and by complexification on $(V_\rr)^\C=E^\C$ that squares to $-$Id. Hence, $E^\C$ splits into the invariant eigen spaces w.r.t. $\mathrm{i} $ can and  $-\i$, which are given by $V$ and $\ol{V}$, respectively. Then
\[x\ni E\ \mapsto\  x-i Jx \ \mapsto \ x+iJ \in \ol{V}\]
gives an isomorphism of real representations
$E\simeq V_\rr\simeq \ol{V}_\rr$.
 Complex multiplication with the imagninary unit on $V$ induces an invariant complex structure $J$ on $E$.

In the other case, where  $E$ is of real type,
$W:=E^\C$ considered as a real vector space,  denoted by $W_\rr$, is reducible,  with invariant real form $E$. This is equivalent to $W$ being {\em self-conjugate} with a conjugation that squares to the identity. Recall that a complex representation $V$ of a real Lie algebra is self conjugate if $W\simeq \ol{W}$ as a $\lag$-representation, i.e. there exists an invariant isomorphism between $W$ and $\ol{W}$. In case of a representation of real type, the invariant real form $E$ then is given as the $+1$ eigen-space of this conjugation.

After this change of  the viewpoint,  it is natural to say  that a {\em complex} irreducible representation is of {\em real type} if it is self-conjugate with a conjugation squaring to one. Otherwise it is called of {\em non-real type}.
Examples of representations of real type are complexifications of the standard representations of $\laso (p,q )$ on $\rr^{p,q}$. Examples of representations of non-real type are representations of $\mf{u}(p,q)$ and $\mf{su}(p,q)$ on $\C^{p,q}$ and $\rr^{2p,2q}$ respectively.

For a complex irreducible representation $V$ of $\lag$ there is a further distinction beyond being of real type or not. If $V$ is not self-conjugate, then $V$ is called {\em of complex type}. If $V$ is self-conjugate with respect to a conjugation $C$, then $C^2$ is a $\C$-linear invariant automorphism of $V$. By the Schur lemma, it is a multiple of the identity, say $C=\lambda \cdot$Id, with $\lambda\in \rr$ because
\[\lambda C (v)= C(\lambda v)=\ol{\lambda }C(v).\]
By scaling $C$ we can assume that $\lambda^2=\pm 1$. In one case, $V$ was of real type, in the case where $C^2=-$Id one says that $V$ is {\em of quaternionic type}, because $C$ defines another complex structure which anti-commutes with the multiplication with i.  To summarise these standard facts, a complex irreducible representation is either of real, complex or quaternionic type. If it not of real type, then $V_\rr$ is irreducible, if it is of real type, it is the complexification of an irreducible real representation.

The following  lemma is a standard result. We cite without proof.
 \blem
$\lag\subset \laso(p,q) $ is of  real type if and only if $p$ and $q$ are even and $\lag\subset \mf{u}(p/2,q/2)$.
\elem
We suppose that the next lemma is also a standard fact from representation theory. Nevertheless, for the sake of being self-contained we prove it here.

\blem\label{oslemma}
Let $\lag $ be a real Lie algebra and $V$ a complex irreducible representation of quaternionic type and of complex dimension $2m$.
\bnum
\item If $V$ is symplectic, then $\lag\subset \mf{sp} (p,q)\subset \mf{u}(2p,2q) $ with $p+q=m$.
\item If $V$ is orthogonal, then $\lag\subset \mf{so}^* (2m)\subset \mf{u}(m,m)$.
\enum
\elem
\bprf
Let $J$ be the anti-linear invariant automorphism of $V$ with $J^2=-1$, and let $V$ be of complex dimension $2m$
Assume that $\w$ is an invariant symplectic form on $V$.
First we show that we can assume the following relation between $\w$ and $J$:
\begin{equation}\w(Jx,Jy)=\ol{\w(x,y)}\label{compatw}\end{equation}
In fact, $\hat{\w}:=\ol{\w(J.,J.)}$ gives another invariant symplectic form on $V$. By the Schur lemma, they are a complex multiple of each other, i.e. $\ol{\w(J.,J.)}=\lambda \w$ for a $\lambda\in \C^*$.  This implies that
\[\w(J.,J.)=\ol{\lambda}\ \ol{\w(J^2.,J^2.)}=\ol{\lambda}\lambda \w(J.,J.)\]
i.e. that $\lambda=\e^{i\theta}\in S^1$. Rescaling $\w$ by $\e^{-i\frac{\theta}{2}}$ enables us to assume equation (\ref{compatw}).
Now note that equation (\ref{compatw}) implies that $\w(J.,.)=-\ol{\w (.,J.)}$. This enables us to define an invariant hermitian form $\la.,.\ra$ on $V$ via
$\la x,y\ra:= \w(x,Jy)$. This is in deed hermitian,
\[\la y,x\ra= \w(y,Jx)=-\ol{\w(Jy,x)}=\ol{\w(x,Jy)}=\ol{\la x,y\ra}\]
and compatible with $J$,
\[\la Jx,Jy\ra = -\w(Jx,y)=\w(y,Jx)=-\ol{\w(Jy,x)}=\ol{\w(x,Jy)}=\ol{\la x,y\ra }.\]
This shows that
 $\lag\subset \mf{u}(2p,2q)\cap \mf{sp}(m,\C)=\mf{sp}(p,q)$, with $p+q=m$.

Now assume that $\s$ is an invariant symmetric bilinear form on $V$. By the same argument as in the symplectic case we get that
\[ \s(Jx,Jy)=\ol{\s(x,y)}\  \text{ and }\  \s(Jx,y)=-\ol{\s(x,Jy)}.\]
A hermitian form is now defined by $\la x,y\ra := \mathrm{i} \s( x,Jy)$.  The compatibility with $J$ is given by
\[\la Jx,Jy\ra = -\mathrm{i}\s(Jx,y) = \mathrm{i}\ol{\s(x,Jy)} = -\ol{\la x,y\ra},\]
which shows that $\la .,.\ra $ has neutral signature $(m,m)$ and that an orthonormal basis of $\s$ is a light-like basis of $\la .,.\ra$.
A calculation in a basis then shows that
 $\lag\subset \mf{u}(V,\la.,.\ra )\cap  \laso (2m,\C )=\laso^* (2m)$, which is defined as follows
 \[\laso^* (2m):=\left\{\left( \begin{array}{cc} A& B\\-\ol{B}^t & \ol{A}^t\end{array}\right)\mid A\in \laso(m,\C), B=\ol{B}^t\right\}\]
(see \cite[p. 446]{helgason78}).
\eprf


\subsection{Irreducibility of  $S^1\cdot SO_0(1,n)$}
In this section we explain that $S^1\cdot SO_0(1,n)\subset U(1,n)\subset  SO(2,2n)$ acts irreducibly.

\bs
 Let $\lag\subset \laso (p,q)$ act irreducibly on $\rr^{p,q}$. This representation is of real type if and only if $\wt{\lag}:=\i\cdot \rr\+ \lag$ acts irreducibly on $\C^{p+q}$ as real representation.
 \es
 \bprf
 If the representation of $\lag$ on $\rr^{p,q}$ is of real type then its complexification is still irreducible, and so is the representation of $\wt{\lag}$ on $\C^{p+q}$. But by definition, there is no conjugation that is invariant under $\wt{\lag}$, and thus
the representation of $\wt{\lag}$ on $\C^{p+q}$ is of complex type, which means that it is still irreducible as a real representation.

On the other hand assume that   $\C^{p+q}$ is irreducible as real and therefore as a complex  representation of $\wt{\lag}$. Assume furthermore that the representation of $\lag$ on $\rr^{p,q}$ is not of real type. By the remarks in the previous section, this is equivalent to the existence of a $\lag$-invariant complex structure $J$ on $\rr^{p,q}$ and to the existence of a complex $\lag$-invariant subspace $V\subset \C^{p+q}$. Extending $J$ complex linear gives an invariant complex structure on $\C^{p+q}$. Note that $J\not= \i\cdot$Id, because otherwise $\wt{\lag}$ could no longer act irreducibly on $\C^{p,q}$. Hence, $I:= \i\cdot J $ is $\wt{\lag}$-invariant, satisfies $I^2=$ Id and is not a multiple of the identity. Hence, it has non-trivial invariant eigen spaces to the eigen values $\pm 1$. But this again contradicts to the irreducibility of $\C^{p,q}$ under $\wt{\lag}$. Therefore, $\rr^{p,q}$ must be of real type for $\lag$.
 \eprf
This gives the following conclusion in the case $p=1$.
 \bfolg
For $n>1$,  $S^1\cdot SO_0(1,n)$ is an irreducible subgroup of $U(1,n)\subset SO_0(2,n)$, is not contained in $SU(1,n)$,  and   has no further irreducible subgroups.
\efolg
\bprf From the previous section we know that irreducible representations of non-real type are unitary, but this is not possible for $\lag \subset \laso(1,n)$. In fact, there is no proper irreducibly acting subalgebra of $\laso(1,n)$, see \cite{olmos-discala01}.  But $\laso(1,n)$ is of real type, and the result follows from the proposition.

For the minimality assume that $\lag\subset\i\rr\+\laso(1,n)$ acts irreducibly. But then the projection of $\lag$ onto $\laso(1,n)$ acts irreducibly and thus has to be equal to $\laso (1,n)$. But this implies $\laso(1,n)=[\lag,\lag]\subset \lag\subset \i\rr\cdot \laso(1,n)$. Hence, $\lag=\i\rr\cdot \laso(1,n)$. That $S^1\times SO(1,n)$ is not contained in $SU(1,n)$ is obvious.
\eprf

\subsection{Reduction to simple Lie algebras}
In this section we show that an irreducible subalgebra of $\laso (2,n)$ is either contained in $\mf{u} (1,n/2)$ or simple.

\bigskip

First, we have to recall some more general facts about representations of real Lie algebras.
Let $\lag$ be a real Lie algebra and $V$ an irreducible {\em complex} representation. We have seen that $V$ is of {\em real type} if and only if there is a $\lag$-invariant existence of a  invariant conjugation squaring to one.
Furthermore, on says that a complex irreducible representation $V$ that is not self-conjugate is of  {\em of complex type}, and if $V$ is self-conjugate with a conjugation squaring to $-1$ it is called {\em of quaternionic type}.
Based on this distinction and on the description of the center in \cite{att05} we proved the following:
 \bs\label{folge3} Let $G\subset SO_0(p,q)$ a
connected Lie subgroup of $SO_0(p,q)$ which acts irreducibly. If $G$
is not semisimple, then $p$ and $q$ are even and $G$ is a subgroup
of $ U(p/2, q/2)$ with centre $U(1)$.
In particular, if $G\subset SO(2,n)$, then $G\subset U(1,n/2)$ or semi-simple.
\es

Here we will strengthen this result for the case $G\subset SO(2,n)$ by replacing ``semi-simple'' by ``simple''. This will be based on the following general fact on complex irreducible representation of semi-simple complex Lie algebras (for a reference, see for example \cite[p. 11]{onishchik04}):
\\[3mm]
{\bf Fact.}
If $\lag=\lag_1\+\lag_2$ is a semi-simple Lie algebra decomposing into non-trivial ideals $\lag_1$ and $\lag_2$.  Then $V$ is a complex irreducible representation of $\lag$ if and only if $V=V_1\otimes V_2$ where $V_i$ are irreducible representations of $\lag_i$.
\\
\blem
Let $\lag\+\lah$ be semi-simple and $W=U\otimes V$ an irreducible complex representation. Then $W$ is self-dual  if and only if both, $U$ and $V$ are  self-dual. The invariant isomorphisms are related by
 $\psi=\psi_1\otimes \psi_2$.
\elem
\bprf
The `if'-direction is obvious, $\psi=\psi_1\otimes \psi_2$ defines the required invariant isomorphism.

For the other direction we consider the identification $\tau: U\simeq U\otimes v_0$ for a fixed $v_0\in V$. $\tau $ is not only an  isomorphism of vector spaces but also of representations of $\lag$, i.e.
\[A (\tau (u))= A (u\otimes v_0)= Au\otimes v_0= \tau (Au).\]
Let $\psi: W\simeq W^*$ be the the isomorphism yielding the self-duality of $W$.
This implies that there are $u_0,\hat{u}_0\in U$ and $v_1\in V$ such that
\[ \left[\psi ( u_0\otimes v_0)\right] (\hat{u}_0\otimes v_1)\not=0.\]
Otherwise, $u_0\otimes v_0  $ would be in the kernel of $\psi$. Hence
by defining
\[\left[ \psi_1 (u)\right](\hat{u}):= \left[\psi ( u\otimes v_0)\right] (\hat{u}\otimes v_1)\]
we obtain a $\lag$-invariant homomorphism $\psi_1:U\simeq U^*$ which is non trivial. By the Schur lemma, $\psi_1$ is an isomorphism.
Obviously, for $V$ one can proceed in the same way.
The Schur-lemma also gives the uniqueness of the invariant structures and the relation between them.
\eprf
\blem
Let $\lag\+\lah$ be semi-simple and $W=U\otimes V$ an irreducible complex self-dual representation. Then $W$ is self-conjugate if and only if both, $U$ and $V$ are  self-conjugate. The invariant isomorphisms are related by
 $\psi=\psi_1\otimes \psi_2$.
\elem
\bprf As $W$ is self-dual,  both $U$ and $V$ are self dual.
Hence, $\ol{U}\simeq \ol{U}^*$ and $\ol{V}\simeq \ol{V}^*$. If $\psi:W\simeq W^*$ and $C:W\simeq \ol{W}$, analogously as in the proof of the previous lemma, one defines  $\phi_1: U\rightarrow \ol{U}^*$ via
\[\left[ \phi_1 (u)\right](\hat{u}):= \left[\psi ( u\otimes v_0)\right] (C(\hat{u}\otimes v_1)).\]
Again, by the Schur lemma, this is an isomorphism, yielding an isomorphism $\psi_1:U\simeq \ol{U}$. All invariant structures are uniquely defined.
\eprf

\btheo
Let $\lag\subset \laso(2,n)$ be an irreducibly acting Lie algebra. Then $\lag\subset \mf{u}(1,n/2)$ or $\lag$ is simple.
\etheo
\bprf
By Proposition \ref{folge3} we can suppose that $\lag$ is semisimple and that the representation of $\lag$ on $\rr^{2,n}$ is of real type.
Assume that
$\lag=\lag_1\+\lag_2$ is not simple. Then its complexification is semisimple and not simple, and thus, the complexified representation $\C^{n+2}$ of  $\rr^{2,n}$ is a tensor product, $\C^{n+2} =V_1\otimes V_2$ of irreducible representations of $\lag_1$ and $\lag_2$. As $\C^{n+2}$ is of real type, the second lemma implies that $V_1$ and $V_2$ are either both of real type or both of quaternionic type. Since $\lag\subset \laso (n+2,\C)$, by the first lemma both are self-dual, defined by either two complex linear symmetric or symplectic forms.

Assume first that both, $V_1$ and $V_2$ are of real type, i.e. $V_i=E_i^\C$ where $E_i$ are irreducible real  representations of $\lag_i$. If $\lag_i\subset \laso (V_i)$, also both $E_i$ are orthogonal, i.e. $\lag_1\subset \laso (p,q)$ and $\lag_2\subset \laso(r,s)$ with
$2=ps+qr$. W.l.o.g. this yields two cases:
The first is $\lag_1=\laso (2)$ and $\lag_2\subset \laso (1,n/2)$ acting on $\rr^2\otimes \rr^{1,\frac{n}{2}}$, or  $\lag_1=\lag_2=\laso(1,1)$. But both cases  contradict to the assumption that $\lag$ was semisimple.

Now we consider the case where the $V_i$'s and thus both $E_i$'s are symplectic representations. In this case the defining scalar product on $\rr^{2,n}$ has neutral signature, i.e. $\lag\subset \laso (2,2)$, and $\lag_i\subset \mf{sp}(1,\rr)=\lasl (2,\rr)$ acting irreducible. Hence, $\lag_i$ either one-dimensional and therefore Abelian, two-dimensional, and thus solvable, or  equal to $\lasl(2,\rr)$. The first two possibilities are excluded by the semisimplicity assumption. We  obtaining that
 $\lag$ is equal to $\lasl(2,\rr)\+\lasl(2,\rr)=\laso (2,2)$.

Now we have to deal with the case where both representations, $V_1$ and $V_2$ are of quaternionic type. As $\lag\subset \laso (n+2,\C)$, they are  either both orthogonal or both symplectic.
Using Lemma \ref{oslemma} we can conclude the proof of the theorem:
First consider the case that $\lag=\lag_1\+\lag_2$ with $\lag_i\subset \mf{sp}(p_i,q_i)\subset \mf{u}(2p_i,2q_i)$. The tensor product of the  hermitian forms on $V_i$ defines a hermitian form of signature $(4(p_1q_2+p_2q_1), 4(p_1p_2+q_1q_2))$ on $V=E^\C$. Since $V$ is an irreducible representation of $\lag$, the space of hermitian forms on $V$ is one-dimensional. Hence, the defined hermitian form is a multiple of the hermitian form obtained by extending the signature $(2,n)$ scalar product on $E$ to $V$. But $2\not= 4(p_1q2+p_2q_1)$ which excludes this case.

For the case $\lag_i\subset \laso^*(2m_i)\subset \mf{u}(m_i,m_i)$ we obtain that $\lag\subset \mf{u}(2m_1m_2, 2m_1m_2)$, which implies $m_i=1$, $V_i=\C^2$ and $\lag_i=\laso^* (2)=\laso(2)$ and $\lag $ is no longer semisimple.
\eprf

 \subsection{Duality for symmetric spaces and consequences}
\label{duality}
For a Riemannian symmetric space $G/K$ given by the Cartan decomposition $\lag=\lak\+\lam$, i.e. by the Lie triple system $\lam$ with the Lie bracket of $\lag$ there is a corresondence between totally geodesic submanifolds and Lie subtriples,
\begin{equation}
\label{Lst}
\{\text{totally geodesic submanifolds $G'/K'$ of $G/K$}\}\simeq
\{\text{Lie subtriples $\lam'\subset \lam$}\}\end{equation}
On the other hand there is the dualtity between compact and non compact symmetric spaces. If $\lag=\lak\+\lam$ is a Cartan decomposition defining  a Riemannian symmetric space, then the one defines the Lie algebra $\lag^*:=\lak\+\lam^*$ by setting
\[\lam^*:= \lam\text{ with Lie bracket }[X,Y]^*:=-[X,Y]\]
and all the other commutators are the same as in $\lag$. In other words, $\lag^*=\lak \+ \i\cdot \lam$ with complex linearly extended commutator.  Then $\lag^*$ defines the dual symmetric space. We get the following correspondence
\be
\{\text{totally geodesic submanifolds  of $G^*/K$}\}&\simeq &\{\text{totally geodesic submanifolds  of $G/K$}\}\\[-2mm]
\turnsimeq\hspace{1cm} &&\hspace{1cm} \turnsimeq\\
\{\text{Lie subtriples in  $\lam^*$}\}&\simeq & \{\text{Lie subtriples in  $\lam$}\}.
\ee
Hence, a totally geodesic submanifold in a compact Riemannian symmetric space $G/K$ is a compact Riemannian symmetric space $H/L$, and the corresponding totally geodescic submanifold in the non-compact dual $G^*/K$ is given by $H^*/L$. For more details we refer the reader to \cite[Chapter 9]{BCO03}.

This correspondence will enable us to describe totally geodesic submanifolds in the Lie Ball $SO_0(2,n)/ SO(2)\times SO(n)$ with the help of totally geodesic submanifolds in the complex quadric $Q^n= SO(n+2)/SO(2)\times SO(n)$.

We can now apply Karpelevich's Theorem \ref{Karp} to what we have obtained so far.
\btheo \label{conseq}
Let $G\subset SO_0(2,n)$ be a connected irreducibly acting subgroup. Then $G\subset U(1,n)$ or $G$ is simple and equal to the effectively acting isometry group of a totally geodesic submanifold in the non-compact symmetric space $SO_0(2,n)/SO(2)\times SO(n)$.
\etheo
\bprf
Let $G\subset SO_0(2,n)$ but $G\not\subset U(1,n)$. From the previous section we know that $G$ is simple. By Karpelvich's Theorem \ref{Karp} it follows that $G$ has a totally geodesic orbit $\cal T$  in the non-compact symmetric space $\cal L^n:= SO_0(2,n)/SO(2)\times SO(n)$. The subgroup
\[ I(\cal T):=\{ A\in G\mid Ap=p\text{ for all } p\in \cal T\}\]
is a normal subgroup in $G$. As $G$ is simple, $I(T)$ is trivial and $G$ acts effectively on $\cal T$. Hence, $\cal T= G/K\subset \cal L^n$ is a non-compact symmetric space with $K\subset G$ maximally compact.
\eprf
In the next section we will determine all irreducibly acting groups $G\subset U(1,n)$ by applying Karpelevich's theorem to the complex projective space. In the last section we will then use a classification of totally geodesic submanifolds in $Q^n= SO(n+2)/SO(2)\times SO(n)$ by \cite{klein08} and the just explained duality to determine the remaining $G$'s.

%

\section{Irreducible subgroups of $U(1,n)$ and complex hyperbolic space}

Using Karpelevich's Theorem in this section we will proof the following statement.

\btheo \label{u1n} Let $G\subset U(1,n)\subset SO(2,2n)$ be a connected subgroup that acts irreducibly on $\rr^{2,2n}$. Then $SU(1,n)\subset G$ or $G=S^1\cdot SO_0(1,n)$.
\etheo
To this end we consider the complex vector space
$\C^{n+1} =: \C^{n,1}$  endowed with the Hermitian form $Q$:

\[ Q =  -|z_0|^2 + |z_1|^2 + |z_2|^2 + \cdots + |z_n|^2 \, . \]

Let us denote by $U(1,n) \subset GL(n+1, \C)$ the subgroup that preserves $Q$.\\

Let $\mathcal{N} := \{ p \in \mathbb{C}^{n,1} : Q(p) < 0 \}$ be the set of negative points. Notice
that $\mathcal{N}$ is a cone preserved by the $U(1,n)$-action. Let us call $\mathbb{C}H^n$ the projectivization of $\mathcal{N}$. Thus, by taking $z_{0} = 1$ we can see that $\mathbb{C}H^n$ is identified with the unit ball of $\mathbb{C}^n$. Namely, \[ \mathbb{C}H^n \cong \{ Z \in \C^n \, : \, |Z|^2 < 1 \} \, \, . \]

It is standard to see that the Hermitian form $Q$ induces on $\ccc H^n$ a $U(1,n)$-invariant Riemannian metric of constant holomorphic curvature. Indeed, we get $\ccc H^n \cong SU(1,n)/U(n)$ as symmetric space of rank one. Notice that the $U(1,n)$-action on $\C H^n$ is not effective since the matrices $e^{i \theta} Id \in U(1,n)$ leaves invariant any complex line. Recall also that the presentation $\ccc H^n \cong SU(1,n)/U(n)$ as symmetric quotient is unique. Namely, if  $\ccc H^n \cong G/K$ where $G$ is semisimple (connected and simply connected) and $K \subset G$ maximal compact then $G = SU(1,n)$ and $K = U(n)$.

The following fact about totally geodesic submanifolds of $\ccc H^n$ can be found in  \cite[pp. 74]{goldman99}, for example.
\bs \label{totalgeo} Let $\mathcal{T} \subset \ccc H^n$ be a complete totally geodesic submanifold. Then $\cal T$ is either a totally real submanifold or a complex submanifold. In the totally real case $\cal T$ is isometric to real hyperbolic space, otherwise $\cal T$ is biholomorphic and isometric to a lower dimensional complex hyperbolic space. In particular,
there exists a real vector subspace $V \subset \ccc^{n,1}$ such that $ \mathcal{T} = V \bigcap \ccc H^n $.
\es

Now we are ready to deduce Theorem \ref{u1n} from Karpelevich's Theorem.

\bprf[Proof of Theorem \ref{u1n}]
Let $H\subset U(1,n)$ be connected and acting irreducibly on $\rr^{2,2n}$ then
 $H$ is reductive, i.e.
 $H=Z\cdot S$ where $Z$ is the centre and $S$ semisimple.
 According with Proposition \ref{folge3} we know that the centre
 $Z$ is trivial or equal to $S^1$.
Hence, the semisimple part $S$ cannot be trivial.

Now,  according to Karpelevich's Theorem $S$ has a totally geodesic orbit $T$ of $\C H^n$. If $T$ is a complex submanifold then Proposition \ref{totalgeo} implies that $S$ must be transitive on $\C H^n$ since otherwise the complex subspace $V$ associated to $T$ is invariant by $S$ and  $Z=S^1$. Thus $H$ can not be irreducible. So $S$ is transitive and we get by the uniqueness of the representation of the symmetric quotient that $SU(1,n) = S$.

Assume now that $T$ is not a complex submanifold. Then, the classification of totally geodesic submanifolds of $\ccc H^n$ imply $T \cong  \R H^n $. Otherwise $T$ is contained in a proper complex totally geodesic submanifold of $\C H^n$ and this imply that $H$ is not irreducible as above. Thus, $T$ is a totally real totally geodesic submanifold.
Without lost of generality we can assume that $T = \R H^n $ where $\R H^n = \{ Z \in \R^n \subset \C^n \, : \, |Z|^2 < 1 \} $. Notice that the Lie algebra of the group $I(T)$ is trivial.
Indeed, if $u \in Lie(I(T))$ then the tangent space to $T$ at $0 \in \R^n \subset \C^n$ is contained in the kernel of $u$. Since $T$ is totally real and $u \in \mf{u}(1,n)$ we get that also the normal space of $T$ at $0$ is contained in the kernel of $u$. Thus $u$ vanish. Since $I(T)$ is trivial we get that $S= SO_0(1,n) \subset SU(1,n)$. Now the center must be $S^1$ and so $H = S^1 \cdot SO_0(1,n) $.
\eprf

\section{The Lie ball and its totally geodesic submanifolds}

 \subsection{The projective model of the Lie ball.}

Let $\R^{2,n}$ the real vector space $\R^{n+2}$ endowed with the
quadratic form
\[q(X,Y) := \langle X, Y \rangle := -x_{0}y_{0} - x_{1}y_{1}+ \sum_{j=2}^{n+1}
x_jy_j, \] where $X = (x_0, \cdots,
x_{n+1})$ and $Y = (y_0, \cdots, y_{n+1})$.
Let $\Pi \subset \R^{2,n}$ be a $2$-dimensional subspace. The
$2$-plane $\Pi$ is called {\it negative} \rm if $q|_{\pi}$ is negative definite.
Let us define the Lie ball ${\cal L}^n$ as one connected component of the set of oriented negative
definite $2$-planes of $\R^{2,n}$. For more details about this model see \cite[p. 285, \S 6]{satake80} or \cite[p. 347]{wolf72}.
Note that $SO(2,n)$ acts transitively on the 
oriented negative
definite $2$-planes, and that $SO_0(2,n)$ acts transitively on
${\cal L}^n$.

Let $\C^{2,n}$ be the complexification of the $\R^{2,n}$, i.e. $q$
becomes \[ q(Z,W) = -z_{0}\overline {w_{0}} - z_{1}\overline{w_{1}} + \sum_{j=2}^{n+1} z_j\overline{w_j}
 ,\] where $Z = (z_0, \cdots, z_{n+1})$ and $W = (w_0, \cdots,w_{n+1})$.
Let $\Pi = span_{\R}\{ A,B\} \subset \R^{2,n}$, $A,B \in
\R^{2,n}$, be an oriented negative definite $2$-plane. We can assume that
$\langle A, B \rangle = 0$ and $q(A,A) = q(B,B) < 0$.
Put $Z = A + iB \, \in \C^{2,n}$ . Then it is not difficult to see
that
\[Z \in Q^{2,n} := \{ Z=(z_0, \cdots , z_{n+1}) \in \C^{2,n} \mid
- z_{0}^2 - z_{1}^2 + \sum_{j=2}^{n+1}z_j^2  = 0 \}\] and that $q(Z,Z) < 0$.
Call $Q^{2,n}_+$ the subset of $Q^{2,n}$ of negative points, i.e.
\[Q_+^{2,n}=\{ Z=(z_0, \cdots , z_{n+1}) \in \C^{2,n} \mid
 - z_{0}^2 - z_{1}^2+\sum_{j=2}^{n+1}z_j^2 = 0\text{ and }q(Z,Z) < 0\}.\]
It follows that we can identify the Lie ball ${\cal
 L}^n$ with a subset of the projective space $\C P^{n,1}$, namely, with a connected component of the image 
of the canonical projection
 $\pi : \C^{n+2} \setminus 0 \rightarrow \C P^{n,1}$.
Thus, we have homogeneous coordinates\footnote{In the appendix an explicit bi-holomorphism with the classical Cartan's domain of type $IV$ is given.} $[z_0 : z_1 : \cdots : z_{n+1} ]$ to work with the Lie ball ${\cal L}^n$.

Let $\Pi_0 = span_{\R} \{ e_{0}, e_{1} \}$ be the ``canonical'' negative definite $2$-plane.
 From now on we will assume that the Lie ball $\cal L^n$
is the connected component of $\Pi_0$.
Then $\Pi_0$ corresponds to the point $Z_0 = e_0 + \i e_{1} = (1,\i, 0,\ldots,0)$.
Thus $\Pi_0 \cong [1: \i: 0:\ldots:0] \in \pi (Q_+^{2,n}) \cong {\cal L}^n $. The isotropy group at $\Pi_0$ is $  SO(2) \times SO(n) $.\\

\subsection{Lifting submanifolds, full submanifolds, and irreducible actions}

Let $M \subset {\cal L}^n$ be a subset. We will denote $L(M) \subset \R^{2,n}$ the subset defined as follows:
\[ L(M) := \bigcup_{ \Pi \in M} \{ p \in \Pi \} \, \, \]
We call $L(M)$ the {\it lift} \rm of $M$.
Let $M \subset \R^n$ be a submanifold. The $M$ is said to be {\it full} \rm if $M$ is not contained in a proper affine subspace of $\R^n$. A submanifold $M \subset {\cal L}^n$ of the Lie ball is called {\it full} \rm if its lift $L(M)$ is full in $\R^{2,n}$.
The following is a well-known property.
\bs Let $G \subset GL(\R,n)$ be a connected Lie subgroup. Assume that $G$ acts irreducibly on $\R^n$.
Then any $G$-orbit $G.p$, $p \neq 0$, is a full submanifold of $\R^n$.
\es

\bprf Let $S_p = span_{\R} \{ G.p \} \subset \R^{n}$ be the linear span of a $G$-orbit $G.p$, $p \neq 0$.
If $S_p$ is not full from some $p$ then $S_p$ lies inside of hyperplane $\langle \cdot, v \rangle = const.$, for some $v \neq 0$. By taking derivatives we get that the Lie algebra $\frak{g}:=Lie(G)$ leaves the subspace $\W = \{ x : \langle x , v \rangle = 0 \}$ invariant. Since $G$ is connected we get that $\W$ is $G$-invariant. \eprf

The following application is also interesting.

\bs \label{full} Let $ G \subset SO(2,n)$ be a connected Lie subgroup and let $\Pi \in {\cal L}^n$ be a point.
Let $G. \Pi$ be the orbit of $G$ through $\Pi$ in the Lie ball.
If the lift $L(G.\Pi)$ is not full then $G$ do not acts irreducibly on $\R^{2,n}$.
\es

\bprf Just notice that the $G$-orbit of any point $p \in \Pi$ is contained in the same proper affine subspace that the lift $L(G.\Pi)$. Then apply the above proposition. \eprf

In the following we will classify full totally geodesic submanifolds in the Lie ball.
This is (almost) equivalent to the classification of {\it maximal} \rm totally geodesic submanifolds in the classical sense (see Sebastian Klein's table at page 11 of \cite{klein08}). We will then check whether the corresponding isometry groups are in our list of irreducible subgroups of $SO_0(2,n)$, respectively, weather or not they are simple (which is the remaining possibility after the previous sections. But first we have to recall the classification of totally geodesic submanifolds in the compact situation, i.e. for the complex quadric $Q^n$.

 \subsection{The complex quadric and its totally geodesic submanifolds}

 The complex quadric $Q^n=SO(n+2)/SO(2)\times SO(n)$ can be viewed in to ways. First, as the Grassmannian of Z(oriented $2$-planes in $\rr^{n+2}$. Secondly, taking into account its complex nature, one can view it as a complex hypersurface in complex projective space, namely as
\[Q^n:=\left\{ [z_0:\ldots : z_{n+1}]\in \C P^{n+1} \mid \sum_{k=0}^{n+1}(z_k)^2=0\right\}.\]
 The subgroup of $SU(n+2)$ acting on $\C^{n+2}$ and thus on $\C P^{n+1}$ that leaves invariant $Q^n$ is $SO(n+2)$ with isotropy group $SO(2)\times SO(n)$. The correspondence to the Grassmannian is given by
\[ P=\span (x,y)\ \mapsto\  \pi (x+iy) \in \C P^{n+1}\]
where $\pi: \C^{n+2}\rightarrow \C P^{n+1}$ is the canonical projection.

Now we will list the totally geodesic submanifolds  in $Q^n$ and their isometry groups as classified in \cite{chen-nagano77} and \cite[Theorem 4.1 and Section 5]{klein08}. Apart from geodesics, there are the following types:

\begin{description}
\item[(I1,\boldmath $k$)] for $1\le k\le n/2 $: This orbit is defined by the following totally geodesic isometric embedding
\[ \C P^k \ni [z_0:\ldots :z_{k}]\ \mapsto\   [z_0:\ldots :z_{k}: \i  z_0:\ldots : \i z_{k}:0:\ldots :0]\in Q^n.\]
Its image is a {\em maximal} totally geodesic submanifolds if $2k=n$ and $n\ge 4$. Its isometry group is $SU(k+1)$ and the totally geodesic submanifold is isometric to $SU(k+1)/U(k)$.
\item[(I2,\boldmath $k$)] for $1\le k\le n/2 $: Here the embedding is give by the restriction of the map for type (I1,k) to real projective space $\rr P^k$ in $\C P^k$. Hence, it is never maximal. Nevertheless, it will be interesting  for our purposes. It is isometric to $O(k+1)/O(k)$.

\item[(G1,\boldmath $k$)] for $1\le k\le n-1$:
This is the embedding of a lower dimensional quadric $Q^k$ into $Q^n$,
\[ Q^k \ni [z_0:\ldots :z_{k+1}]\ \mapsto\   [z_0:\ldots :z_{k+1}:0:\ldots :0]\in Q^n.\]
It is maximal for $k=n-1\ge 2$. Its isometry group is $SO(k+2)$ and it is isometric to $SO(k+2)/SO(2)\times SO(k)$.
\item[(G2,\boldmath $k_1$,$k_2$)] for $ 1\le k_1+k_2 \le n$: This is a totally geodesic isometric embedding of a product of two spheres with radius $1/\sqrt{2}$ and of dimension $k_1$ and $k_2$ given by
\[
\left((x_0,\ldots , x_{k_1}) ,  (y_0,\ldots , y_{k_2})\right) \mapsto
 \left[ x_0:\ldots :x_{k_1}: \i  y_0:\ldots : \i y_{k_2}:0:\ldots :0\right] \in Q^n
 \]
 This orbit is maximal for $k_1+k_2=n\ge 3$. Its isometry group is given by
 $SO(k_1+1)\times SO(k_2+1)$.
 \item[(G3)] The quadric $Q^2$ is isometric to $\C P^1\times \C P^1$ i.e., $\C P^1\times \C P^1 \equiv Q^2  $ .
 Let $C=\mathbb{R}P^1 \subset \C P^1$ be the trace of a closed geodesic in $\C P^1$.
 Then the map

 \[ \C P^1\times C \rightarrow  \C P^1\times \C P^1 \equiv Q^2 \rightarrow Q^n \, \, \]

 where the last embedding represents the embedding of type (G1,2) described above. So the embedding $\C P^1\times C \hookrightarrow Q^m $ is maximal only for $n=2$.

  \item[(P1,\boldmath $k$)] for $1\le k\le n$. This is given as the embedding of type (G2, $k_1$, $k_2$) for $k_1$ or $k_2$ equal to zero. Its image is maximal for $k=n$. The isometry group is given as $SO(k+1)$.
  \item[(P2)] This is the embedding of type (G1, $k$) for $k=1$. It is maximal only for $n=2$ and its isometry group is $SO(3)$.
 \item[(A)] The totally geodesic submanifold is isometric to the $2$-sphere of radius $\sqrt{10}/2 $. It is maximal only for $n=3$ and its isometry group is given by $SO(3)$.
\end{description}

\subsection{Totally geodesic submanifolds of the Lie ball} \label{type}

We will now use Cartan's duality (as explained in Section \ref{duality}) and Klein's classification as listed in the previous section.
In the following, the immersions $u$ will be equivariant.
So they are useful to compute the corresponding immersion of the group into $SO(2,n)$.

\paragraph{Type (I1,k)} Here we have $1\le k\le n/2$.
Let us consider the following map,
 \[ u\ :\ [z_0 : \ldots : z_k] \longrightarrow [z_0 : \i z_0 : \ldots  : z_k : \i z_k:0:\ldots :0]. \]
The image of $u$ is contained in $\pi(Q^{2,n})$. In order to see which point is taken by $u$  to ${\cal L}^n $ it is enough to see that
\[ -|z_0|^2 - |\i z_0|^2 + \sum_{i=1}^{n+1}\left( |z_i|^2 + |\i z_i|^2\right)  = 2  \left(- |z_0|^2 + \sum_{i=1}^{n+1} |z_i|^2 \right) \]
Thus, $u[z_0 :\ldots :z_k:0:\ldots :0] \in {\cal L}^n $ if and only if $-|z_0|^2 + \sum_{i=1}^{n+1} |z_i|^2 < 0$.
Hence,  $u$ gives an holomorphic immersion from the complex hyperbolic space $\C H^k$ into our Lie ball ${\cal L}^n$. Namely, $\C H^k$ is regarded as the projective submanifold of $\C P^{k,1}$ with $-|z_0|^2 + \sum_{i=1}^{n+1} |z_i|^2 < 0$\footnote{To see $\C H^k$ as a ball in $\C^k$ just take $z_0 = 1$, i.e. the affine chart.}.

The group of isometries of $\C H^k$ is $SU(1,k) \subset SO(2,n)$ which acts irreducibly on $\R^{2,n}$ only for $k=n/2$.
To see this it is enough to identify $\R^{2,n}$ with $\C^{1,n/2}$ endowed with the quadratic form $-|w_0|^2 + \sum_{i=1}^{n/2}|w_i|^2$. The action of $SU(1,n/2)$ is transitive on the set of negative $2$-planes of $\C^{1,n/2}$ given by complex lines. For example, the complex line generated by the vector $(1,0,\ldots , 0) \in \C^{1,n/2}$ is a negative definite $2$-plane\footnote{Actually, such a complex line is the $2$-plane which we called $\Pi_0$ in the first section.} of $\R^{2,n}$.
Let $w = (w_0, \ldots ,w_{n/2}) \in \C^{1,n/2}$ be a vector. Then the $2$-plane generated by $w$, i.e. the complex line, is given by the homogeneous coordinates $[ \overline{w_0} : \i \overline{w_0} : \ldots  : \overline{w}_{n/2} : \i \overline{w}_{n/2} ]$. This show that the image of our map $u$ is the set of $2$-planes coming from complex lines of $\C^{1,n/2}$. Thus, the image $u(\C H^{n/2})$ is the orbit of $SU(1,n/2)$ through $\Pi_0$.

Notice that the lift $L(u(\C H^{n/2}))$ of the totally geodesic submanifold $u(\C H^{n/2})$ is
the union of the points in all negative complex lines. Thus, such a subset is full in $\R^{2,n}$ and
this is consistent (indeed equivalent) to the fact that $SU(1,n/2)$ acts irreducibly on $\R^{2,n}$.

\paragraph{Type  (I2,\boldmath $k$)} Here it is  $1\le k\le n/2 $.
The map $u$ is the ``real form'' of the above map:

 \[ [x_0 : \ldots  : x_k] \stackrel{u}{\rightarrow} [x_0 : \i x_0 :\ldots  : x_k : \i x_k:0:\ldots :0] \]

Thus we get an embedding  of $\R H^k$ (in the projective Klein model\footnote{Here we refer to Felix Klein.}) into the Lie ball.
Notice that the subgroup $SO(1,k) \subset SU(1,n/2)$ acts reducibly on $\R^{2,n}$, even for $k=n/2$. In the light of Theorem \ref{conseq} we do not get another irreducible subgroup of $SO(2,n)$.
But we should point out that
the group $I(u(\R H^k))$ i.e. the isometries that fix all points of the image of $u$
is given as  $I(u(\R H^k)) = SO(2)$ acting diagonally, i.e. $SO(2) \cong e^{\operatorname{i} \theta} Id$. For $k=n/2$  this group makes $G=I(u(\R H^k))\cdot SO(1,n/2)$ act irreducibly on $\rr^{2,n}$. $G$ was already on our list.

\paragraph{Type (G1,\boldmath $k$)} This is the embedding of a lower dimensional Lie ball. Its isometry group is given by $SO(2,k)$, which does not act irreducibly on $\rr^{2,n}$.

\paragraph{Type (G2,\boldmath $k_1,k_2$)} In this case  $1\le k_1+k_2 \le n$ and the map $u$ is given by:
\[ \left( [x_0:\ldots : x_{k_1}],[y_0:\ldots :y_{k_2}]\right) \stackrel{u}{\mapsto} [x_0 : \i y_0 : x_1:\ldots : x_{k_1} :\i y_1 :\ldots : \i y_{k_2} :0 : \ldots :0] \]
The image lies in ${\cal L}^n$ if and only if:
\[ -x_0^2 - y_0^2 +\sum_{i=1}^{k_1}x_i^2 + \sum_{j=1}^{k_2}y_j^2 < 0 \]
and
\[- x_0^2 + y_0^2 +\sum_{i=1}^{k_1}x_i^2 - \sum_{j=1}^{k_2}y_j^2=0. \]
Since the map is given in homogeneous coordinates we can assume that
$-x_0^2 +\sum_{i=1}^{k_1}x_i^2 = - y_0^2 + \sum_{j=1}^{k_2}y_j^2$ which shows that
the image of $u$ is in the Lie ball if and only if
$[x_0:\ldots : x_{k_1}]$ and $[y_0:\ldots :y_{k_2}]$ lie in the real hyperbolic spaces of dimensions $k_1$ and $k_2$. Hence, $u$ is an embedding of $\R H^{k_1}\times \R H^{k_2}$ into the Lie ball. The isometry groups is given by $SO(1,k_1)\times SO(1,k_2) \subset SO(2,k_1+k_2)\subset SO(2,n)$.  Thus, the isometry group of this totally geodesic submanifold does not act irreducibly, since it fixes $\R^{1,k_1}$ and $\R^{1,k_2}$.

\paragraph{Type (P1,\boldmath $k$)} Here it it $1\le k\le n$ and the embedding is given by the one of type (G2, $k_1$, $k_2$) for $k_1$ or $k_2$ equal to zero. Hence, we can write it as
\[[ x_0:\ldots : x_k ] \rightarrow [\i: x_0:\ldots : x_k:0\ldots :0] \]
Thus we get an immersion from $\R H^k$ (as in the usual Lorentzian model) into the Lie ball ${\cal L}^n$. The isometry group of the totally geodesic submanifold is $SO(1,k)$ acting reducibly even for $k=n$ by fixing the first basis vector $e_0$.

\paragraph{Type (P2)} This is the embedding of type (G1, $k$) for $k=1$ and thus  the isometry group of the totally geodesic submanifolds is given as $SO(2,1)\subset SO(2,n)$ acting reducibly by fixing $e_3, \ldots, e_{n+1}$.

\paragraph{Type (G3)} This totally geodesic submanifold is a Riemannian product. Then its isometry group $G$ is not simple. Thus this case reduce to the case of  $G \subset U(1,n)$.

\paragraph{Type (A)} Here it is $n\ge 3$. This is an embedding of $3$-dimensional real hyperbolic space into the Lie ball. The only irreducible acting (simple) subgroup of $SO(2,n)$ which did not  appear as isometry group of a totally geodesic orbit in the Lie ball is $SO(1,2)\subset SO(2,3)$. Thus we conclude that this embedding of $SO(1,2) $ gives the isometry group of a totally geodesic orbit of type $(A)$ for $n=3$. For $n>3$ it is reducible, of course. For details on this case please refer to Appendix \ref{so12}.

\bigskip

We conclude that the only irreducibly acting simple proper subgroups of $SO(2,n)$ that appear as isometry group of a totally geodesic submanifold in the Lie ball are the following
\[  SU(1,n/2) \text{ and } SO(1,2)\subset SO(2,3).\]

\section{Proof of  Theorem \ref{list}}

Let $G \subset SO(2,n)$ be a connected subgroup whose action on $\R^{2,n}$ is irreducible. Assume that $G \neq SO_0(2,n)$.
If $G$ is not simple then Theorems \ref{conseq} and \ref{u1n} imply that $G$ is one of the groups in our list. Namely, $n=2p$ and either $G = U(1,p)$ or $G=S^1\cdot SO_0(1,p)$ if $p>1$. If $G$ is simple then (up to conjugation) $G$ is the group of isometries of one of the totally geodesic submanifolds of the Lie ball ${\cal L}^n$ listed in the previous section. Thus, either $G = SU(1,n/2)$ or $ G = SO_0(1,2)\subset SO(2,3)$. This completes the proof.\hfill $\Box$

\begin{appendix}

\section{Appendix}

\subsection{\boldmath $\laso(1,2)$ acting irreducibly on \boldmath $\rr^{2,3}$ and orbits of type (A)}
\label{so12}
In this section we want to describe the irreducible injections of $\laso (3)$ and $\laso (1,2)$ into $\laso (5)$ and $\laso (2,3)$, respectively, and to describe the Lie subtriples in the complex quadric and the Lie ball corresponding to them.

The irreducible injection of $\laso (3)\subset \laso (5)$ corresponds to an irreducible symmetric space  of type $AI$ for $n=3$, which we will describe in the non-compact setting. To this end split $\lasl (3,\rr)$ as
\[\lasl (3,\rr)=\laso (3)\+ sym_0(3)\]
where $sym_0(3)$ denotes the trace free, symmetric $3\times 3$ matrices. This splitting provides symmetric data, since
\[ [\laso(3),\laso(3)]\subset \laso(3),\  [\laso(3),sym_0(3)]\subset sym_0(3),\ [sym_0(3),sym_0(3)]\subset \laso(3).\]
The adjoint representation of $\laso (3) $ on the five-dimensional space $sym_0(3)$ is orthogonal with respect to the Killing form of $\lasl(3,\rr)$, which is positive definite on $sym_0(3)$. Hence, the above splitting defines an irreducible $5$-dimensional Riemannian symmetric space with isotropy group $\laso (3)\subset \laso(5)$, irreducibly.

In order to write down $\laso (3)$ in $5\times 5$ matrices,
we identify the standard basis in $\rr^5$ with the following matrices which are orthogonal under the Killing form of $\lasl (3,\rr)$:
\[e_1= S_{12},\ e_2= S_{13},\ e_3=\mathrm{diag}(-1,1,0),\  e_4=S_{23},\ e_5=\frac{1}{\sqrt{3}}\ \mathrm{ diag}(-1,-1,2),\]
where $S_{ij}$ denotes the symmetric matrix with $1$ at the $(i,j)$-th spot.
Acting via the adjoint representation, the standard basis of $\laso (3)$
\[ U=D_{12}\ ,\ \ V=D_{13}\ ,\text{ and } W=D_{23}\]
is given as follows
\[
U=
\left(\begin{array}{c|c} 0&\begin{array}{ccc} -2&0&0\\0&-1&0\end{array}\\ \hline
 \begin{array}{cc}2&0\\
 0&1\\
 0&0\end{array} &0
 \end{array}\right) \text{, }
V=
\left(\begin{array}{c|c} 0&\begin{array}{ccc} 0&-1&0\\-1&0&-\sqrt{3}\end{array}\\ \hline
 \begin{array}{cc}0&1\\
 1&0\\
 0&\sqrt{3}\end{array} &0
 \end{array}\right)
\]
and
\[
W\ =\  [U,V]\ =\
\left(\begin{array}{c|c}
\begin{array}{cc}
0&-1\\
1&0
\end{array}
&0 \\ \hline
0&
 \begin{array}{ccc}0&-1&0\\
 1&0&-\sqrt{3}\\
 0&\sqrt{3}&0\end{array}
 \end{array}\right)
\]
One verifies that
\be
[V,W]= U&\text{ and } &[W,U]=V,\ee
which are the commutator relations of $\laso (3)$.

Regarding the complex quadric $Q^n=SO(n+2)/SO(n)$,
the Cartan decomposition of the Lie algebra $\laso (n+2)$ is given by
 \be
\lak\ :=\  \laso(2)\+\laso(n)&=&\left\{\left(\begin{array}{c|c} A&0\\ \hline 0&B\end{array}\right)\mid A\in \laso(2), B\in \laso(n)\right\}\\
 \lam&:=&\left\{\left(\begin{array}{c|c} 0&\begin{array}{c} -u^t\\-v^t\end{array}\\ \hline
 \begin{array}{cc}u&v\end{array} &0
 \end{array}\right)\mid u,v\in \rrn\right\}
 \ee
Now, for $n\ge 3$, $\lan =\span (U,V)$ defines a Lie sub-triple of the Lie triple corresponding to the complex quadric, with isometry algebra $\laso (3)= [\lan,\lan]\+\lan\subset \laso (5)\subset \laso (n+2)$. The totally geodesic orbit corresponding to this Lie subtriple is of type (A), see \cite[page 85]{klein08}.

\bigskip

Now we transfer this situation to the Lie ball with Cartan decomposition
\[\laso(2,n)=\lak \+\lam^*\ \text{ with }\
 \lam^*:=\left\{\left(\begin{array}{c|c} 0&\begin{array}{c} u^t\\v^t\end{array}\\ \hline
 \begin{array}{cc}u&v\end{array} &0
 \end{array}\right)\mid u,v\in \rrn\right\}
\]
Here the totally geodesic orbits of type (A) are given by the Lie subtriple
\[\lan^*= \span (U^*,V^*)\]
where $U^*$ and $V^*$ are as above, but symmetric instead of skew symmetric. Here we have that
\[ -W=[U^*,V^*]\]
 and again
\[ [V^*,W]= U^*\ \text{ and }\  [W,U^*]=V^*.\]
These are the commutator relations of $\laso(1,2)\simeq \lasl(2,\rr)$. Again, the irreducible representation $\laso(1,2)\subset \laso (2,3)$ comes from the irreducible symmetric space of signature $(2,3)$ given by
\[\lasl (3,\rr)=\laso(1,2)\+ \mathfrak{t} \]
where $\mathfrak{t} $ is a five-dimensional complement of $\laso(1,2) $ in $\lasl(3,\rr)$, consisting of  trace-free matrices with the right symmetries.

Finally, in order to verify a remark in the introduction, we want to show that $SO_0(1,2)$ does not act locally transitively on $S^{1,2}=SO_0(2,3)/P$  where $P$ is the parabolic subgroup given as the stabiliser of a light-like line. For general $(p,q)$, the Lie algebra of $P$ is given as $\lap=(\rr\+\laso(p,q))\ltimes \rr^{p,q}\subset \laso(p+1,q+1)$. A group $G\subset SO(p+1,q+1)$ acts locally transitively on $S^{p,q}=SO(p+1,q+1)/P$ if
\[\laso (p+1,q+1)=\lag + \lap\]
where $\lag\subset \laso(p+1,q+1)$ is the Lie algebra of $G$ (for details see \cite{alt08}). For $\lag=\laso (1,2)\subset \laso(2,3)$, irreducibly, this sum has to be direct since the parabolic $\lap$ is $7$-dimensional in this case. But $U^*\in \laso(1,2) $ fixes the line spanned by the light-like vector vector $e_1+e_3$. Hence, $\laso(1,2)\cap\lap \not=\{0\}$ and thus, this action of $SO(1,2)$  on $S^{1,2}$ is not  locally transitively.

\subsection{A biholomorphism between the Lie ball and  Cartan's bounded domain of type $IV$}
In this appendix we give an explicit bi-holomorphism between the Lie ball ${\cal L}^n$ and (the classical) Cartan's bounded domain of type $IV$ in $\C^n$.
Let $f$ be the map given as follows:

\[ f(z_1, \cdots, z_n) = [i(\Lambda - 1) : \Lambda + 1 : 2z_1 : \cdots, 2z_n  ] \, \, ,\]

where $\Lambda := z_1^2 + \cdots + z_n^2$. Thus, $f: \C^n \rightarrow \C P^{1,n}$. Let us show that $f(\C^n) \subset Q^{2,n}$. Indeed,  \[ (2z_1)^2 + \cdots + (2z_n)^2 - (i(\Lambda - 1))^2 - (\Lambda + 1)^2 = \]
\[ 4 \Lambda + (\Lambda -1)^2 - (\Lambda + 1)^2 = 4 \Lambda + (\Lambda - 1 - \Lambda - 1)(\Lambda - 1 + \Lambda + 1) = \] \[ = 4 \Lambda - 4 \Lambda = 0 \]

Let us characterize when a point $z = (z_1, \cdots, z_n)$ is take by $f$ to a point in the Lie ball ${\cal L}$.
By definition $f(z) \in {\cal L}^n$ if and only if the following condition holds:

\[ |2z_1|^2 + \cdots + |2z_n|^2 - |i(\Lambda - 1)|^2 - |\Lambda + 1|^2 < 0 \]

Notice that $f(0) = [ -i : 1: 0: \cdots : 0 ] = [1 : i: 0: \cdots : 0 ] = \Pi_0$ satisfy such a condition. Thus, we are interested in the connected component of $0 \in \C^n$ satisfying such a condition. A simple calculation shows that the above condition is equivalent with:

\[ 2(|z_1|^2 + \cdots + |z_n|^2) - |z_1^2 + \cdots + z_n^2|^2 - 1 < 0 \]

Notice that $|z_1^2 + \cdots + z_n^2|^2 \leq (|z_1|^2 + \cdots + |z_n|^2)^2$. So let $z(t)$ be a continuous curve starting at $0$ (i.e., $z(0)=0$) and satisfying the above condition. Then at each $t$ we get
\[ 2(|z_1|^2 + \cdots + |z_n|^2)  - 1 <  |z_1^2 + \cdots + z_n^2|^2 \leq (|z_1|^2 + \cdots + |z_n|^2)^2 \]
So
\[ 2(|z_1|^2 + \cdots + |z_n|^2)  - 1 <  (|z_1|^2 + \cdots + |z_n|^2)^2 \]
\[ - (|z_1|^2 + \cdots + |z_n|^2  - 1)^2 < 0 \]

Thus the curve $z(t)$ remains inside the unit ball $|z| < 1$.\\

This shows that $f$ takes the bounded domain $\Omega \subset \C^n$ given by the inequalities:

\[ |z| < 1 \]

and

\[ 2(|z_1|^2 + \cdots + |z_n|^2) - |z_1^2 + \cdots + z_n^2|^2 - 1 < 0 \]

into the Lie ball ${\cal L}^n$. Notice that $\Omega$ is indeed the classical Cartan's bounded domain of type $IV$.
Actually, $f$ is just the stereographic projection of the quadric $- w_{1}^2 - w_{0}^2 + w_1^2 + \cdots + w_n^2  = 0$ regarded as an affine sphere in the chart $w_{0} = 1$ of $\C P^{1,n}$.
\end{appendix}


\def\cprime{$'$}

\end{document}